\input amstex.tex 
\documentstyle{amsppt}
\pagewidth{6.5truein}
\pageheight{8.6truein}

\topmatter
\title Small Zeros of Quadratic Forms with Linear Conditions \endtitle
\author Lenny Fukshansky \endauthor
\address Department of Mathematics, The University of Texas at Austin, Austin, Texas 78712 \endaddress
\email lenny\@math.utexas.edu \endemail
\subjclass Primary 11D09, 11E12; Secondary 11H46 \endsubjclass
\abstract Given a quadratic form and $M$ linear forms in $N+1$ variables with coefficients in a number field $K$, suppose that there exists a point in $K^{N+1}$ at which the quadratic form vanishes and all the linear forms do not. Then we show that there exists a point like this of relatively small height. This generalizes a result of D.W. Masser.
\endabstract
\endtopmatter

\document
\flushpar

\loadbold

\def\M{{\Cal M}}

\def\s{{\Cal S}}

\def\V{{\Cal V}}

\def\be{{\boldsymbol e}}

\def\bq{{\boldsymbol q}}
\def\bs{{\boldsymbol s}}
\def\bt{{\boldsymbol t}}
\def\bu{{\boldsymbol u}}

\def\bw{{\boldsymbol w}}
\def\bx{{\boldsymbol x}}
\def\bX{{\boldkey X}}
\def\bz{{\boldsymbol z}}
\def\bwy{{\boldsymbol y}}
\def\bY{{\boldkey Y}}

\def\ba{{\boldsymbol\alpha}}

{\bf \S1. Introduction and notation.} Let
$$F(\bX,\bY) = \sum_{i=0}^{N} \sum_{j=0}^{N} f_{ij} X_i Y_j$$
be a symmetric bilinear form in $N+1$ variables with coefficients $f_{ij} = f_{ji}$. We write $F = (f_{ij})$ for the associated $(N+1) \times (N+1)$ matrix, and $F(\bX) = F(\bX, \bX)$ for the associated quadratic form. First assume that the coefficients $f_{ij}$ are in $\Bbb Q$. Suppose there exists a point $\bx \in \Bbb Q^{N+1}$ such that $x_0 \neq 0$ and $F(\bx) = 0$. In \cite{4} Masser shows that in this case there exists such a point $\bx$ with
$$H(\bx) \ll_N H(F)^{(N+1)/2},$$
where $H$ here stands for height of $\bx$ and $F$, respectively. This generalizes a well known result of Cassels \cite{2} about the existence of small zeros of quadratic forms with rational coefficients to the existence of small zeros of quadratic polynomials with rational coefficients.

We generalize Masser's result in the following way. Let $K$ be a number field of degree $d$ over $\Bbb Q$. Let the coefficients $f_{ij}$ be in $K$. Let $M$ be a positive integer. Let $L_1(\bX),...,L_M(\bX)$ be linear forms in $N+1$ variables with coefficients in $K$. Suppose there exists a point $\bt \in K^{N+1}$ such that $F(\bt) = 0$, and $L_i(\bt) \neq 0$ for each $1 \leq i \leq M$. Then we prove that there exists such a point of bounded height. The bound on height is in terms of the heights of quadratic and linear forms, and reduces (up to a constant) to Masser's type result over a number field in case $M=1$ and $L_1(\bX)=X_0$. 
\bigskip

First we set some notation. For a number field $K$ of degree $d$ over $\Bbb Q$, we write $O_K$ for the ring of algebraic integers of $K$ and $\Delta_K$ for the discriminant of $K$. Write $M(K)$ for the set of all places of $K$, and for each $v \in M(K)$ let $d_v = [K_v:\Bbb Q_v]$ be the local degree, where $K_v$ and $\Bbb Q_v$ are completions of $K$ and $\Bbb Q$ respectively at the place $v$. Then if $u \in M(\Bbb Q)$, let $M_u = \left\{ v \in M(K) : v|u \right\}$, and we have
$$\sum_{v \in M_u} d_v = d.$$
We normalize our absolute values for $v \in M(K)$ as in \cite{6}:
\roster
\item if $v | p$ then $|p|_v = p^{-d_v/d}$,
\item if $v | \infty$ then $|\alpha|_v = |\alpha|^{d_v/d}$, where $|\ |$
is the usual Euclidean absolute value on $\Bbb R$ or $\Bbb C$.
\endroster

\noindent
Then for every $\alpha \in K$, $\alpha \neq 0$, the product formula reads
$$\prod_{v} |\alpha|_v = 1.$$

\noindent
For each $v \in M(K)$, we define a local height over $K_v$ by
$$H_v(\bx) = \max_{0 \leq i \leq N} |x_i|_v,$$
for each $\bx \in K_v^{N+1}$. Then we have the {\it homogeneous} global height function on $K^{N+1}$
$$H(\bx) = \prod_{v \in M(K)} H_v(\bx),$$
and the {\it inhomogeneous height}
$$h(\bx) = \prod_{v \in M(K)} \max\{1,H_v(\bx)\},$$
for each $\bx \in K^{N+1}$. We now state a basic well-known property of height functions. Let $\bx,\bwy \in K^N$, and $\alpha,\beta$ be positive integers, then
$$H(\alpha\bx \pm \beta\bwy) \leq h(\alpha\bx \pm \beta\bwy) \leq (\alpha+\beta)h(\bx)h(\bwy).\tag1.1$$
We define the height of a polynomial to be the height of its coefficient vector.
\smallskip

\noindent
Let $j$ be a positive integer. Let
\roster
\item $r_v(j) = \pi^{-1/2} \Gamma(j/2+1)^{1/j}$, if $v | \infty$ is real, 
\item $r_v(j) = (2\pi)^{-1/2} \Gamma(j+1)^{1/2j}$, if  $v | \infty$ is complex,
\item $r_v(j)=1,$ if $v \nmid \infty$.
\endroster
For each $j$, define a field constant
$$A_K(j) = \left\{ 2^{5j} (j+1)^j |\Delta_K|^{\frac{j+1}{d}} \right\}^{1/2} \prod_{v \in M(K)} r_v(j)^{\frac{jd_v}{d}}.\tag1.2$$
\smallskip

Now we can rigorously state the main result of this paper.

\proclaim{Theorem 1.1} Suppose there exists a point $\bt \in K^{N+1}$ such that $F(\bt) = 0$, and $L_i(\bt) \neq 0$ for each $1 \leq i \leq M$. Then there exists $\bu \in K^{N+1}$ such that $F(\bu) = 0$, $L_i(\bu) \neq 0$ for each $1 \leq i \leq M$, and
$$H(\bu) \leq B_K(N,M) H(F)^{\frac{N+2M}{2}+(M-1)(N+2)},\tag1.3$$
as well as
$$H(\bu) \leq B_K(N,M) H(F)^{\frac{N+1}{2}+(M-1)(N+2)} \prod_{i=1}^M H(L_i)^{\frac{(2M-1)N}{M}},\tag1.4$$
and finally
$$H(\bu) \leq B_K(N,M) H(F)^{\frac{2N+2M+1}{4}+(M-1)(N+2)} \prod_{i=1}^M H(L_i)^{\frac{(2M-1)N}{2M}},\tag1.5$$
where the constant $B_K(N,M)$ is given by
$$B_K(N,M) = \frac{1}{192} (N+1)^2 A_K(N) \left\{ 486\ (N+1)^6 A_K(N)^2 \right\}^{M-1} (M+2)! \{(M+3)!\}^2,\tag1.6$$
with $A_K(N)$ as in (1.2).
\endproclaim

The following result is a simple, but useful corollary of Theorem 1.1 in the case $M=1$.

\proclaim{Corollary 1.2} Let $F(\bX)$ be a quadratic form in $N+1$ variables with coefficients in the number field $K$, as above. Let
$$\V_K(F) = \left\{\bt \in K^{N+1} : F(\bt)=0 \right\}.$$
Suppose that there exists a non-singular point $\boldkey 0 \neq \bx \in \V_K(F)$. Then there exists a non-singular point $\boldkey 0 \neq \bs \in \V_K(F)$ such that
$$H(\bs) \leq \max \{3, A_K(N)\}\ H(F)^{\frac{N}{2}}.$$
\endproclaim

The structure of this paper is the following. In \S2 we produce a solution to the problem in case there is only one linear form, obtaining upper bounds for the inhomogeneous height of the point in question, and proving Corollary 1.2. Our line of argument here follows that of Masser \cite{4}. In the process of proof we state a generalization of Cassels' result on small zeros of quadratic forms, that we use to construct auxiliary points. In \S3 we produce an upper bound for the height of a point outside of the collection of subspaces. In \S4 we prove Theorem 1.1. It is derived from a slightly more technical result of Theorem 4.1. Our argument is by induction on the number of linear forms, so we use the results of \S2 for the base case of the induction, and we use the result of \S3 to construct certain auxiliary points. Then we compute bounds on the height. We also remark that one can assume the point $\bu$ of Theorem 1.1 to be in $O_K^{N+1}$.
\bigskip

{\bf \S2. The problem with one linear form.} Let $L(\bX)$ be a linear form in $N+1$ variables with coefficients in $K$, and suppose there exists a point $\bt \in K^{N+1}$ so that $F(\bt) = 0$ and $L(\bt) \neq 0$. We want to show an existence of such a point of small height. The argument of this section parallels that of Masser \cite{4}. We argue by induction on $N$.

First suppose that $N=1$, then
$$F(X_0,X_1) = aX_0^2 + bX_0X_1 + cX_1^2,$$
$$L(X_0,X_1) = q_0X_0 + q_1X_1,$$
where $a,b,c,q_0,q_1 \in K$. Since $L(\bX)$ is not identically zero, we can assume without loss of generality that $q_0 \neq 0$. If $a=c=0$, then $\alpha \{(1,0),(0,1)\}$, $\alpha \in K$, is the zero set of $F$ consisting of two projective points of height $1$, and $L$ must not vanish at one of them. Then assume $a \neq 0$. Let $\bx = (x_0,x_1) \in K^2$ be a non-trivial zero of $F$, so $x_0,x_1 \neq 0.$ Then
$$(x_0,x_1) = x_1 \left(\frac{-b \pm \sqrt{b^2-4ac}}{2a}, 1\right),$$
so again the zero set of $F$ consists of only two projective points, and hence $L$ must not vanish at one them. Thus we just have to estimate the heights of these two points. We can assume that $x_1=1$, therefore $h(\bx)=H(\bx)$. A straightforward calculation shows that
$$h(\bx) \leq 3 H(F),\tag2.1$$
and this finishes the proof in case $N=1$.
\bigskip

Now we state a generalized form of Cassels' theorem on small zeros of quadratic forms, that we will use in the proof. The following version is due to Vaaler.

\proclaim{Theorem 2.1} If a quadratic form $F$ has a nontrivial zero in $K^{N+1}$, then there exists $\boldkey 0 \neq \bx \in O_K^{N+1}$ such that $F(\bx) = 0$, and
$$H(\bx) \leq h(\bx) \leq A_K(N) H(F)^{N/2},\tag2.2$$
where $A_K(N)$ is as in (1.2).
\endproclaim

This follows by combining Theorem 1, Corollary 2 and remark after it of \cite{6} with Corollary 11 of \cite{1}.
\smallskip

{\it Remark.} A theorem like this has first been proved for the case $K=\Bbb Q$ by Cassels in \cite{2}, and later generalized to number fields by Raghavan \cite{5} (various other important generalizations of Cassels' result were also carried out by Birch, Davenport, Chalk, Schmidt, Schlickewei, and Vaaler, just to name a few; see \cite{6} for a more detailed account and bibliography).
\bigskip

We return to the proof. Now assume that $N \geq 2$. Then 
$$L(\bX) = \bq \centerdot \bX = \sum_{i=0}^N q_iX_i \in K[X_0,...,X_N].\tag2.3$$
By Theorem 2.1, there exists $\boldkey 0 \neq \bx \in K^{N+1}$ such that $F(\bx) = 0$ and
$$H(\bx) \leq h(\bx) \leq A_K(N) H(F)^{N/2}.\tag2.4$$
If $L(\bx) \neq 0$, we are done, so assume $L(\bx) = 0$. Again, since $L(\bX)$ is not identically zero, we can assume that for instance $q_0 \neq 0$.  
This implies that
$$x_0 = -\frac{1}{q_0} \sum_{i=1}^{N} q_ix_i,$$
hence
$$0 = F(\bx) = \sum_{i=1}^N \sum_{j=1}^N f_{ij}x_ix_j + 2 \sum_{i=1}^N f_{0i}x_0x_i + f_{00}x_0^2 = \sum_{i=1}^N \sum_{j=1}^N g_{ij} x_ix_j,$$
where for each $1 \leq i,j \leq N$, $g_{ij} = f_{ij} - \frac{2q_j}{q_0}f_{0i} + \frac{f_{00}}{q_0^2}q_iq_j$.
Then define a quadratic form $G$ in $N$ variables:
$$G(\bX) = \sum_{i=1}^{N} \sum_{j=1}^{N} g_{ij} X_i X_j.$$
Notice that $\boldkey 0 \neq (x_1,...,x_N) \in K^N$, and $G(x_1,...,x_N) = 0$, hence by Theorem 2.1, there exists $\boldkey 0 \neq \bz \in K^N$ such that $G(\bz) = 0$ and
$$H(\bz) \leq h(\bz) \leq A_K(N-1) H(G)^{(N-1)/2}.$$
We need a bound on $H(G)$ in terms of $H(F)$ and $H(L)$. Using the fact that $H_v(L) \geq |q_0|_v$ for each $v \in M(K)$ along with ultrametric inequality in the non-archimedean case and triangle inequality in the archimedean case, we obtain
$$H(G) \leq 6 H(F) H(L)^2,$$ 
therefore
$$H(\bz) \leq h(\bz) \leq 6 A_K(N-1) H(F)^{(N-1)/2} H(L)^{N-1}.\tag2.5$$
Define
$$y_0 = -\frac{1}{q_0} \sum_{i=1}^{N} q_iz_i,$$
and let $\boldkey 0 \neq \bwy = (y_0, \bz) \in K^{N+1}$. By construction, $F(\bwy) = L(\bwy) = 0$. Then using (2.5), we obtain
$$H(\bwy) \leq h(\bwy) \leq N \prod_{v \in M(K)} \frac{H_v(L)}{|q_0|_v} \max \{1,H_v(\bz)\}\tag2.6$$
$$\leq 6 N A_K(N-1) H(F)^{(N-1)/2} H(L)^N.$$
\smallskip

Since the bilinear form $F$ is not identically zero, there must exist a coefficient $f_{ij} \neq 0$. Then without loss of generality, assume $f_{00} = 1$, which implies that
$$\max \{1, H_v(F)\} = H_v(F),\tag2.7$$
for each $v \in M(K)$.
\smallskip

Next let $\boldkey 0 \neq \bt_1, \bt_2 \in K^{N+1}$, and define
$$\bu_1 = F(\bt_1)\bx - 2F(\bt_1,\bx)\bt_1,\tag2.8$$
and
$$\bu_2 = F(\bt_2)\bwy - 2F(\bt_2,\bx)\bt_2.\tag2.9$$
It is easy to check that $F(\bu_1) = F(\bu_2) = 0$. Let
$$\V_K(F) = \left\{ \bt \in K^{N+1} : F(\bt) = 0 \right\}.$$

\proclaim{Lemma 2.2} Suppose that $\bx, \bwy$ are non-singular points in the variety $\V_K(F)$. Then there exist $\boldkey 0 \neq \bt_1, \bt_2 \in K^{N+1}$ with coordinates $0, \pm 1$ such that $L(\bu_1), L(\bu_2) \neq 0$.
\endproclaim

\demo{Proof} We will go through the construction of $\bt_1$, and the construction of $\bt_2$ is identical. Since $L(\bx) = 0$, we want to construct $\bt_1 \in K^{N+1}$ such that the following holds:
\roster
\item $t_{10} \neq -\frac{1}{q_0} \sum_{i=1}^N q_it_{1i},$
\item $F(\bt_1,\bx) \neq 0,$
\item $t_{1i} = 0,\pm 1\ \ \forall\ \  0 \leq i \leq N.$
\endroster
Notice that (1) is equivalent to $L(\bt_1) \neq 0$, and (2) is possible since $\bx$ is non-singular in $\V_K(F)$. Write $\be_0,...,\be_N$ for the standard basis vectors. Each $\be_i$ satisfies (3). There exists $\be_i$ satisfying (1). If $\be_i$ satisfies (2), let $\bt_1 = \be_i$. Otherwise, there exists $\be_j$ satisfying (2), and $i \neq j$. If $\be_j$ satisfies (1), let $\bt_1 = \be_j$. If not, then let $\bt_1 = \be_i + \be_j$, and we are done.  
\boxed{ }
\enddemo
\bigskip

Assume $\bx, \bwy$ are non-singular points in the variety $\V_K(F)$. Make the choice of $\bt_1, \bt_2$ in (2.8), (2.9) as in Lemma 2.2. Then $F(\bu_1) = F(\bu_2) = 0$, $L(\bu_1), L(\bu_2) \neq 0$. We want to estimate heights of $\bu_1,\bu_2$.

\proclaim{Lemma 2.3} If $\bt, \bw \in K^{N+1}$, and $\bu = F(\bt)\bw - 2F(\bt,\bw)\bt$, then
$$H(\bu) \leq h(\bu) \leq 3(N+1)^2H(F)h(\bw)h(\bt)^2.\tag2.10$$
\endproclaim

\demo{Proof} If $v \nmid \infty$, then $|2|_v \leq 1$, and so
$$\max \{1,H_v(\bu)\} \leq \max \left\{1, |F(\bt)|_vH_v(\bw), |2|_v|F(\bt,\bw)|_vH_v(\bt)\right\}$$
$$\leq \max \{1, H_v(F)H_v(\bw)H_v(\bt)^2\} \leq \max \{1, H_v(F)\} \max \{1, H_v(\bw)\} \max \{1, H_v(\bt)\}^2$$
$$= H_v(F) \max \{1, H_v(\bw)\} \max \{1, H_v(\bt)\}^2,$$
where the last equality follows by (2.7). If $v | \infty$, then
$$H_v(\bu) \leq |F(\bt)|_vH_v(\bw) + 2|F(\bt,\bw)|_vH_v(\bt) \leq \{3(N+1)^2\}^{d_v/d} H_v(F)H_v(\bw)H_v(\bt)^2,$$
and so
$$\max \{1,H_v(\bu)\} \leq \{3(N+1)^2\}^{d_v/d} \max \{1, H_v(F)H_v(\bw)H_v(\bt)^2\}$$
$$\leq \{3(N+1)^2\}^{d_v/d} \max \{1, H_v(F)\} \max \{1, H_v(\bw)\} \max \{1, H_v(\bt)\}^2$$
$$= \{3(N+1)^2\}^{d_v/d} H_v(F) \max \{1, H_v(\bw)\} \max \{1, H_v(\bt)\}^2,$$
where the last equality follows by (2.7). Then (2.10) follows by taking a product.
\boxed{ }
\enddemo
\smallskip

By Lemma 2.2, $h(\bt_1) = h(\bt_2) = 1$, and so by Lemma 2.3, (2.4), and (2.6) we have
$$h(\bu_1) \leq 3(N+1)^2 H(F)h(\bx) \leq 3 (N+1)^2 A_K(N) H(F)^{(N+2)/2},\tag2.11$$
and
$$h(\bu_2) \leq 3(N+1)^2 H(F)h(\bwy) \leq 18 N(N+1)^2 A_K(N-1) H(F)^{(N+1)/2}H(L)^N.\tag2.12$$
\bigskip

Next we consider the ``singular'' case.

\proclaim{Proposition 2.4} Assume that $\bx$ is a singular point in the variety $\V_K(F)$. Then there exists a point $\bs \in K^{N+1}$ so that $F(\bs) = 0$, $L(\bs) \neq 0$, and
$$H(\bs) \leq h(\bs) \leq 3 H(F)^{N/2}.\tag2.13$$
\endproclaim

\demo{Proof} Here the idea is as in \cite{4}, to reduce to fewer variables keeping coefficients under control and to use induction. If $N=1$, (2.13) is just (2.1). Then assume that $N \geq 2$, and that (2.13) has been proved for $N-1$. Without loss of generality, assume that $x_N \neq 0$. Then $\bx$ is linearly independent of the first $N$ standard unit vectors $\be_0,...,\be_{N-1}$, so we can define new variables $Y_0,...,Y_N$ by
$$\bX = (X_0,...,X_N) = Y_0\be_0+...+Y_{N-1}\be_{N-1}+Y_N\bx.\tag2.14$$
We have
$$F(\bX) = F\left( \sum_{i=0}^{N-1} Y_i\be_i \right) + Y_N^2 F(\bx) + 2F\left( \sum_{i=0}^{N-1} Y_i\be_i, Y_N\bx \right) = F\left( \sum_{i=0}^{N-1} Y_i\be_i\right),$$
since $F(\bx)=0$, and $\bx$ is a singular point in $V$, i.e. $F(\bt,\bx)=0$ for all $\bt \in K^{N+1}$. Then define a new quadratic form $Q$ in $N$ variables $Y_0,...,Y_{N-1}$ by
$$Q(\bY) = F\left( \sum_{i=0}^{N-1} Y_i\be_i \right),$$
and so $F(\bX)=Q(\bY)$. Clearly, the coefficients of $Q$ form a subset of coefficients of $F$, and hence
$$H(Q) \leq H(F).\tag2.15$$
There exists a $\bt \in K^{N+1}$ so that $F(\bt) = 0$, and $L(\bt) \neq 0$. Let $\bw = (w_0,...,w_{N-1})$ be the vector that corresponds to $\bt$ under the coordinate change (2.14) and reduction to $N$ variables. Then
$$0 \neq L(\bt) = L \left(\sum_{i=0}^{N-1} w_i\be_i \right) + \frac{t_N}{x_N} L(\bx) = L \left(\sum_{i=0}^{N-1} w_i\be_i \right),$$
since $L(\bx)=0$. Then define a new linear form $L_1$ in $N$ variables $Y_0,...,Y_{N-1}$ by
$$L_1(\bY) = L \left(\sum_{i=0}^{N-1} Y_i\be_i \right),$$
and so $L_1(\bw) \neq 0$, and
$$H(L_1) \leq H(L),$$
since coefficents of $L_1$ form a subset of coefficients of $L$. We also know that $Q(\bw) = F(\bt) = 0$. Therefore, by induction hypothesis, there exists $\bu \in K^N$ such that $Q(\bu)=0$, $L_1(\bu) \neq 0$, and
$$h(\bu) \leq 3 H(Q)^{N/2} \leq 3 H(F)^{N/2},$$
by (2.15). Define $\bs = (\bu,0) \in K^{N+1}$, and then $F(\bs)=Q(\bu)=0$, $L(\bs)=L_1(\bu) \neq 0$, and $h(\bs) = h(\bu)$. This completes the proof.   
\boxed{ }
\enddemo
\smallskip

Now notice that if $N \geq 1$, $3(N+1)^2 A_K(N) > 3$, as well as for each $N \geq 2$, $(N+1)^2 A_K(N) \geq N(N+1)^2 A_K(N-1)$. Putting this together with (2.11), (2.12), and (2.13), we have proved the following theorem.

\proclaim{Theorem 2.5} Let the notation be as above. Suppose there exists a point $\bt \in K^{N+1}$ such that $F(\bt) = 0$, and $L(\bt) \neq 0$. Then there exists $\bu \in K^{N+1}$ such that $F(\bu) = 0$, $L(\bu) \neq 0$, and
$$H(\bu) \leq h(\bu) \leq 18(N+1)^2 A_K(N) H(F)^{(N+1)/2} \min \left\{H(F)^{1/2}, H(L)^N \right\}.\tag2.16$$
\endproclaim
\smallskip

\demo{Proof of Corollary 1.2} Let $\bx$ be the zero of $F$ guranteed by Theorem 2.1. If $\bx$ is non-singular, we are done. If $\bx$ is singular, let $L(\bX) = \frac{\partial F}{\partial X_i}$ for some $0 \leq i \leq N$, so $L(\bx) = 0$. Then by Proposition 2.4, there must exist $\bs \in K^{N+1}$ so that $F(\bs) = 0$, $L(\bs) \neq 0$, and
$$H(\bs) \leq h(\bs) \leq 3 H(F)^{N/2}.$$
\boxed{ }
\enddemo
\bigskip

{\bf \S3. Points of small height outside of a collection of subspaces.} Let $M,N$ be positive integers, and let $K$ be a number field of degree $d$ over $\Bbb Q$. Keeping all the notation as before, we prove the existence of a point of small height at which none of $M$ linear forms in $N$ variables with coefficients in $K$ vanish.
\smallskip

In fact, we consider a more general situation and produce a basic result. Let $v \in M(K)$ be any place of $K$, and let 
$$\s_M = \{ \ba \in \Bbb Z_{\geq 0}^N : \alpha_1+...+\alpha_N \leq M \}.$$
Then let
$$U(X_1,...,X_N) = \sum_{\ba \in \s_M} c(\ba) X_1^{\alpha_1}...X_N^{\alpha_N} \in K_v[X_1,...,X_N],$$
be a polynomial in $N$ variables of degree $M$. If $k$ is a positive integer, then for each vector $\bx \in \Bbb Z^k$ write
$$|\bx| = \max \{|x_1|,...,|x_k|\}.$$
The idea for the following argument was suggested to me by Sinnou David, \cite{3}.

\proclaim{Theorem 3.1} Let the notation be as above, and suppose $U(\bX)$ is not identically $0$. Then there exists $\bx \in \Bbb Z^N$ such that $U(\bx) \neq 0$ and $|\bx| \leq \left[\frac{\deg(U)}{2}\right]+1 = \left[\frac{M}{2}\right]+1$.
\endproclaim

\demo{Proof} We argue by induction on $N$. First suppose $N=1$. Then our polynomial is of the form
$$U(X) = c_M X^M+...+c_1 X+c_0 \in K_v[X],$$
and $U$ has at most $M$ integer roots. Hence there must exist $x \in \Bbb Z$ such that $U(x) \neq 0$ and $|x| \leq \left[\frac{M}{2}\right]+1$. Now suppose the theorem has been proved for all polynomials in $k$ variables for any $1 \leq k < N$. Notice that for each $1 \leq i \leq N$, $\deg_{X_i}(U) \leq \deg(U)=M$, where $\deg_{X_i}(U)$ is the degree of $U$ in the variable $X_i$.

There must exist $\bq \in \Bbb Z^{N-1}$ such that $U(\bq, X_N)$ is not identically $0$. Indeed, suppose it is not so. Then $U$ vanishes on all of $\Bbb Z^N$, which by continuity implies that $U$ is identically $0$.

Since $U(\bq,X_N)$ is a polynomial in one variable, by the base of induction there exists $q_N \in \Bbb Z$ such that $U(\bq, q_N) \neq 0$ and $|q_N| \leq \left[\frac{M}{2}\right]+1$. Let
$$P(X_1,...,X_{N-1}) = U(X_1,...,X_{N-1},q_N),$$
then $P$ is not identically $0$, and $\deg(P) \leq M$. By induction hypothesis, there exists $\bx \in \Bbb Z^{N-1}$ such that $P(\bx) \neq 0$ and $|\bx| \leq \left[\frac{M}{2}\right]+1$. Then $(\bx,q_N) \in \Bbb Z^N$, $U(\bx,q_N) = P(\bx) \neq 0$, and $|(\bx,q_N)| \leq \left[\frac{M}{2}\right]+1$.
\boxed{ }
\enddemo
\bigskip

\proclaim{Corollary 3.2} Let the notation and assumptions be as in Theorem 3.1. Then there exists $\bx \in O_K^N$ such that $U(\bx) \neq 0$ and $H(\bx) \leq h(\bx) \leq \left[\frac{\deg(U)}{2}\right]+1 = \left[\frac{M}{2}\right]+1$.
\endproclaim

\demo{Proof} The point $\bx$ obtained in Theorem 3.1 is in $\Bbb Z^N \subseteq O_K^N$, and so
$$h(\bx) \leq \prod_{v | \infty} \max\{1,|\bx|_v\} = \prod_{v | \infty} \max\{1,|\bx|\}^{d_v/d} = |\bx|.$$
The result follows.
\boxed{ }
\enddemo
\smallskip

Considering the special case when $U(\bX) = \prod_{i=1}^M L_i(\bX)$ is decomposable into a product of $M$ linear forms $L_1,...,L_M$, we conclude that there exists a point $\bx \in K^N$ at which none of the linear forms vanish and $h(\bx) \leq \frac{M+2}{2}$.
\bigskip

{\bf \S4. Proof of Theorem 1.1.} Let $M$ and $N$ be positive integers. Let $F$ be a quadratic form in $N+1$ variables with coefficients in a number field $K$ of degree $d$, as above. Let $L_1,...,L_M$ be linear forms in $N+1$ variables with coefficients in $K$.

\proclaim{Theorem 4.1} Suppose there exists a point $\bt \in K^{N+1}$ such that $F(\bt) = 0$, and $L_i(\bt) \neq 0$ for each $1 \leq i \leq M$. Then there exists $\bu \in K^{N+1}$ such that $F(\bu) = 0$, $L_i(\bu) \neq 0$ for each $1 \leq i \leq M$, and
$$H(\bu) \leq h(\bu) \leq B_K(N,M) H(F)^{\frac{N+1}{2}+(M-1)(N+2)} \prod_{i=1}^M \M_i^{2-\frac{1}{M}},$$
where
$$\M_i = \min \left\{H(F)^{1/2}, H(L_i)^N \right\},\tag M$$
and the constant $B_K(N,M)$ is as in (1.6).
\endproclaim

\demo{Proof} We will actually prove a slightly stronger upper bound:
$$h(\bu) \leq B_K(N,M) H(F)^{\frac{N+1}{2}+(M-1)(N+2)} \M_1 \prod_{i=2}^M \M_i^2.\tag*$$
We argue by induction on $M$. If $M=1$, then Theorem 4.1 follows from Theorem 2.5. So suppose $M \geq 2$, and that theorem has been proved for any subset of $L_1,...,L_M$ of $k$ linear forms, where $1 \leq k \leq M-1$. Then there exist points $\bx,\bwy \in K^{N+1}$ such that $F(\bx) = F(\bwy) = 0$, $L_i(\bx) \neq 0$ for every $1 \leq i \leq M-1$, $L_M(\bwy) \neq 0$, and
$$h(\bx) \leq B_K(N,M-1) H(F)^{\frac{N+1}{2}+(M-2)(N+2)} \M_1 \prod_{i=2}^{M-1} \M_i^2,\tag4.1$$
$$h(\bwy) \leq 18(N+1)^2 A_K(N) H(F)^{(N+1)/2} \M_M.\tag4.2$$
Notice that if $b < a$ are positive integers, we interpret $\prod_{i=a}^b$ as $1$. If $L_M(\bx) \neq 0$ or $L_i(\bwy) \neq 0$ for all $1 \leq i \leq M-1$, then we are done. So assume it is not so. Then there exists a $k$, such that $1 \leq k < M-1$ and by reordering the linear forms if necessary we have
\roster
\item $L_i(\bx) \neq 0$, $L_i(\bwy) \neq 0$, for all $1 \leq i \leq k$,
\item $L_i(\bx) \neq 0$, $L_i(\bwy) = 0$, for all $k < i \leq M-1$,
\item $L_M(\bx) = 0$, $L_M(\bwy) \neq 0$.
\endroster 
Notice that for every $k < i \leq M$, $L_i(\bx+\bwy) \neq 0$. In fact, there exists a positive integer $\beta$ such that for all $1 \leq i \leq M$,
$$L_i(\bx \pm \beta\bwy) \neq 0,$$
for the same choice of $\pm$. For this, $\beta$ needs to be such that for the same choice of $\pm$ none of the linear equations in $\beta$
$$L_i(\bx) \pm \beta L_i(\bwy) = 0,\ \ 1 \leq i \leq k \leq M-2,$$
are true. There are at most $M-2$ such equations, and since we can also choose $\pm$, there exists such a $\beta$ so that
$$1 \leq \beta \leq \left[ \frac{M-2}{2} \right]+1 \leq \frac{M}{2}.\tag4.3$$
Define
$$\bu = \bx \pm \beta\bwy,$$
for this choice of $\pm$ and $\beta$.
\smallskip

{\it Case 1.} Suppose $F(\bx,\bwy)=0$. Then
$$F(\bu)=F(\bx)+\beta^2 F(\bwy) \pm 2\beta F(\bx,\bwy) = 0,$$
and
$$L_i(\bu) \neq 0,\ \ \forall\ 1 \leq i \leq M.$$
Combining (1.1) and (4.3) we obtain
$$h(\bu) \leq (\beta+1)h(\bx)h(\bwy) \leq \left( \frac{M+2}{2} \right) h(\bwy) h(\bx).\tag4.4$$
\bigskip

{\it Case 2.} Suppose $F(\bx,\bwy) \neq 0$. By Corollary 3.2, there exists $\bw \in K^{N+1}$ such that $L_i(\bw) \neq 0$ for each $1 \leq i \leq M$ and
$$h(\bw) \leq \frac{M+2}{2}.\tag4.5$$
If $F(\bw) = 0$, we are done. Assume it is not so. Let $\beta$ be a positive integer, and define
$$\bu = F(\bwy \pm \beta \bw)\bx - 2F(\bx, \bwy \pm \beta \bw)(\bwy \pm \beta \bw).$$
Notice that $F(\bu) = 0$. We want to choose $\pm \beta$ in such a way that the following is true:
\roster
\item $F(\bwy \pm \beta \bw) = \beta (\beta F(\bw) \pm 2F(\bwy,\bw)) \neq 0$,
\item $F(\bx, \bwy \pm \beta \bw) = F(\bx,\bwy) \pm \beta F(\bx,\bw) \neq 0$,
\item $L_i(\bu) = F(\bwy \pm \beta \bw)L_i(\bx) - 2F(\bx, \bwy \pm \beta \bw)(L_i(\bwy) \pm \beta L_i(\bw)) \neq 0$, for each $1 \leq i \leq M$.
\endroster
It is not difficult to see that (1), (2), (3) amount to a total of $2$ linear and $M$ quadratic expressions in $\beta$. Selecting $\pm$ appropriately we see that there exists a positive integer $\beta$ such that (1), (2), (3) are satisfied, and
$$\beta \leq M+2.\tag4.6$$
By the same argument as in \S2, we can assume without loss of generality that $f_{00} = 1$. Then, for this choice of $\pm \beta$, Lemma 2.3, (1.1), (4.5), and (4.6) imply that
$$h(\bu) \leq 3(N+1)^2 H(F)h(\bx)h(\bwy \pm \beta \bw)^2\tag4.7$$
$$\leq 3(N+1)^2 (\beta + 1)^2 \left( \frac{M+2}{2} \right) H(F)h(\bx)h(\bwy)^2$$
$$\leq \frac{3}{2}(N+1)^2 (M+2)(M+3)^2 H(F)h(\bx)h(\bwy)^2.$$
Combining (4.2), (4.4), and (4.7), we have proved that there exists $\bu \in K^{N+1}$ such that $F(\bu)=0$, $L_i(\bu) \neq 0$ for each $1 \leq i \leq M$, and
$$h(\bu) \leq 486 (N+1)^6 A_K(N)^2 (M+2)(M+3)^2 H(F)^{N+2} \M_M^2 h(\bx).\tag4.8$$
This proves (*). Notice that the ordering of linear forms was arbitrary, so assume that $\M_1 = \max_{1 \leq i \leq M}\ \M_i$. Then $\M_1 \geq \M_1^{1/M}...\M_M^{1/M}$, and so
$$\M_1 \prod_{i=2}^M \M_i^2 \leq \prod_{i=1}^M \M_i^{2-\frac{1}{M}}.\tag4.9$$
The theorem follows by combining (4.8) with (4.1) and (4.9).
\boxed{ }
\enddemo
\smallskip
 
To derive Theorem 1.1, notice that for each $1 \leq i \leq M$, the following inequalities hold:
$$\M_i \leq H(F)^{1/2},\ \ \ \M_i \leq H(L_i)^N,\ \ \ \M_i \leq H(F)^{1/4} H(L_i)^{N/2}.$$
Combining these with the inequality of Theorem 4.1 produces (1.3), (1.4), and (1.5) respectively. 
\bigskip

{\it Remark.} Let $\bu$ be the point of Theorem 1.1. Since $H_v(\bu) = 1$ for all but finitely many places of $K$, the Strong Approximation Theorem guarantees the existence of $0 \neq \alpha \in K$ such that $\alpha \bu \in O_K^{N+1}$, and of course $H(\alpha \bu) = H(\bu)$. From this, it is also an easy exercise to produce an upper bound on $H_{\infty}(\alpha \bu) = \prod_{v | \infty} H_v(\alpha \bu)$. The exponents in the upper bound turn out to be the same as in Theorem 1.1, and the constant is only slightly larger.
\bigskip

{\bf Aknowledgements.} I would like to express my deep gratitude to Professor Jeffrey D. Vaaler for his valuable advice and numerous helpful conversations on the subject of this paper. I would also like to thank Professor Sinnou David for his helpful idea that I used in \S4, and the referee for his useful comments and simplifications of some arguments.


\Refs
\ref \no 1 \by E. Bombieri, J. D. Vaaler \pages 11--32
\paper On Siegel's lemma
\yr1983
\jour Invent. Math.\vol 73
\endref
\ref \no 2 \by J. W. S. Cassels \pages 262--264
\paper Bounds for the least solutions of homogeneous quadratic equations
\yr1955
\jour Proc. Cambridge Philos. Soc. \vol 51
\endref
\ref \no 3 \by S. David
\jour personal communication
\endref
\ref \no 4 \by D. W. Masser \pages 24--28
\paper How to solve a quadratic equation in rationals
\yr1998
\jour Bull. London Math. Soc.\vol 30
\endref
\ref \no 5 \by S. Raghavan \pages 109--114
\paper Bounds of minimal solutions of diophantine equations
\yr1975
\jour Nachr. Akad. Wiss. Gottingen, Math. Phys. Kl.\vol 9
\endref
\ref \no 6 \by J. D. Vaaler \pages 281--296
\paper Small zeros of quadratic forms over number fields
\yr1987
\jour Trans. Amer. Math. Soc.\vol 302
\endref

\endRefs
\enddocument